\documentclass{article}

\usepackage[T2A]{fontenc}
\usepackage[utf8]{inputenc}
\usepackage[english,russian]{babel}
\usepackage[tbtags]{amsmath}
\usepackage{amsfonts,amssymb}
\usepackage{amsbib}

\newcommand{\im}{\operatorname{Im}}
\newcommand{\re}{\operatorname{Re}}

\begin{document}

\author{ В.\,Е.~Владыкина \thanks{МГУ им. М.\,В.~Ломоносова, e-mail: vika-chan@mail.ru} ,
 А.\,А.~Шкаликов\thanks {МГУ им. М.\,В.~Ломоносова, e-mail: shkaliko@yandex.ru}
}
 \title{Асимптотика решений уравнения Штурма--Лиувилля с сингулярными коэффициентами}
\date{6.10.2015\thanks{Работа опубликована в Математических заметках, {\bf 98} (6), 832-841 (2015)}}

\maketitle
{\sc УДК 517.928 + 517.984}
\begin{abstract}
Получены асимптотические представления  при $\lambda \to \infty$
 в верхней и нижней полуплоскостях  для решений уравнения Штурма-Лиувилля
 $$
 -y''+p(x)y'+q(x)y= \lambda ^2 \rho(x)y, \qquad x\in [a,b] \subset  \mathbb{R},
 $$
 при условии, что $q$ --- распределение первого порядка сингулярности, $\rho$ ---
положительная  абсолютно непрерывная функция, а $p$  принадлежит пространству $L_2[a,b]$.

Библиография: 8 названий.
\end{abstract}
\footnotetext{Работа выполнена при финансовой
поддержке фонда РНФ, \No~14-11-00754.}

Представления для асимптотических решений дифференциальных уравнений по спектральному параметру в секторах комплексной плоскости хорошо известны и часто используются начиная с работ Биркгофа \cite {1},\cite{2}. Конечно, впоследствии появились многочисленные работы, где результаты Биркгофа развивались и обобщались. Отметим, в частности, монографии Тамаркина \cite{3}, Рапопорта \cite{4} и Наймарка \cite{5} (в последней монографии имеется достаточно подробная библиография). Однако вплоть до недавнего времени все результаты на эту тему относились только к уравнениям, в которых коэффициенты при младших членах являются классическими функциями (непрерывными или суммируемыми по Лебегу). Нам известна только одна работа, в которой ограничения такого рода не предполагались. Это работа Савчука и Шкаликова \cite{6}, в которой были предприняты различные подходы для определения оператора Штурма--Лиувилля с потенциалом-распределением первого порядка сингулярности. В частности, в \cite{6} были получены асимптотические представления решения уравнения
\begin{equation} \label{1}
-y''+q(x) y =\lambda ^2 y
\end{equation}
при условии, что потенциал $q$ является обобщенной функцией, принадлежащей  пространству $W_2^{-1}[a,b]$, $[a,b] \subset \mathbb{R}$. Последнее означает, что понимаемая в смысле распределений первообразная $\int q \, dt$ принадлежит пространству $L_2[a,b]$. Конечно, полученные асимптотические представления существенно использовались для представления функции Грина соответствующего оператора и изучения его спектральных свойств. Они использовались также и при изучении обратных задач (см., например, \cite{7}).

Напомним в упрощенной формулировке результат \cite{6} об асимптотическом представлении для фундаментальной системы решений уравнения \eqref{1}.

{\sc Теорема} \cite{6}. {\it
При условии $q \in W_2^{-1}[a,b]$ уравнение \eqref{1} имеет пару линейно независимых решений вида
	\begin{equation} \label{2}
	y_s(x, \lambda )= \sin \lambda x + \varphi(x, \lambda), \quad	y_c(x, \lambda )= \cos \lambda x + \psi(x, \lambda),
	\end{equation}
	где $ |\varphi(x, \lambda)|+ |\psi(x, \lambda)| = o(1)$
 при $|\lambda | \to \infty$ в полосе $|\im \lambda| \leqslant r$  равномерно по $x \in [a,b]$. В качестве $r$ можно взять произвольное положительное число. Асимптотические представления \eqref{2}, вообще говоря, нельзя почленно дифференцировать, но если ввести квазипроизводную  $y^{[1]}=y'-u(x)y$, где $u=\int q \, dt$, то
	$$	y^{[1]}_s(x, \lambda )=\lambda( \cos \lambda x + \varphi_1(x, \lambda)), 	
\quad y^{[1]}_c(x, \lambda )= -\lambda ( \sin \lambda x + \psi _1(x, \lambda))
$$
с сохранением свойства $ |\varphi_1(x, \lambda)|+ |\psi _1(x, \lambda)| = o(1)$ при $|\lambda | \to \infty$ в полосе $|\im \lambda| \leqslant r$  равномерно по $x \in [a,b]$.
}

	Упрощение приведенной здесь формулировки  состоит в том, что, в отличие от \cite{6}, мы не детализируем
специфические свойства функций $\varphi$ и $\psi$,  ограничиваясь  лишь их основным свойством:
стремлением к нулю при $|\lambda | \to \infty$ в полосе произвольной фиксированной ширины.
Отметим, что доказательство сформулированной теоремы потребовало в \cite{6} весьма сложной технической работы.

	Цель настоящей работы -- получить асимптотические представления для решений уравнения более общего вида
	\begin{equation}\label{3}
	-y''+p(x)y'+q(x)y= \lambda ^2 \rho(x)y, \qquad x\in [a,b] \subset  \mathbb{R},
	\end{equation}
	при следующих условиях на коэффициенты
\begin{equation} \label{4}	
p(x) \in L_2[a,b], \quad q \in W_2^{-1}[a,b], \quad\rho \in W^1_1[a,b],\ \
\rho(x)>0 \ \ \text{при} \ \, x\in [a,b].
\end{equation}		
Здесь $ W^1_1[a,b]$ -- пространство функций, для которых $y' \in L_1[a,b]$ (в силу теоремы Лебега это пространство совпадает с пространством $AC[a,b]$ абсолютно непрерывных функций). Условие $q \in W_2^{-1}[a,b]$ означает, как уже отмечалось, что первообразная
$u(x) = \int q(t)\, dt$,  понимаемая в смысле теории распределений, принадлежит
пространству $L_2[a,b]$. Функции $p$ и $q$ предполагаются комплекснозначными.
Кроме того, будем предполагать выполненным также условие
 \begin{equation} \label{5}
\rho '(x) u(x) \in L_1[a,b],  \quad \text{где} \  u(x) =  \int q(x) \, dx.
\end{equation}

  Помимо того, что мы в сравнении с \cite{6} рассматриваем более общее уравнение, мы имеем еще и другую цель. Мы намерены получить асимптотические представления для решений не только в полосе, но и в каждой из полуплоскотей  $\im \lambda \geqslant -r$ и
   $\im \lambda \leqslant r$,  где $r\geqslant 0$ произвольное фиксированное число. Это существенное добавление. Наличие асимптотических представлений в областях, которые покрывают всю комплексную плоскость, является фактом, который, несомненно, полезен для дальнейших исследований. Отметим также, что в настоящей работе мы предлагаем другой метод получения асимтпотических представлений, нежели в \cite{6}. Этот метод, по-видимому, имеет более широкие возможности.

Сформулируем и докажем основной результат работы.

{\sc Теорема.} {\it  При выполнении условий \eqref{4} и \eqref{5}
	уравнение \eqref{3} имеет на отрезке $[a,b]$ фундаментальную систему решений вида
\begin{equation} \label{6}
y_{\pm}(x, \lambda)= \left[\rho(x)\right]^{-1/4}\exp\left( \frac 12 P(x) \pm i \lambda\, \int_a^x\sqrt{\rho(\tau)}\, d\tau\right)
\left(1+ \varphi_\pm(x, \lambda)\right).
 	\end{equation}
	Здесь
$$
P(x)=\int\limits_a^x p(\tau)\, d\tau,
$$
 а функции $\varphi_{\pm}$ таковы, что
 $$
  |\varphi _+(x, \lambda)|+ |\varphi _-(x, \lambda)| = o(1) \quad \text{при} \ \, |\lambda | \to \infty \ \,
  \text{в полуплоскости} \ \, \im \lambda \geqslant -r,
  $$
  равномерно по $x\in [a,b]$, причем число $r$ можно взять произвольно большим. Асимптотики \eqref{6} можно почленно дифференцировать, если вместо производной рассматривать квазипроизводную
	 $$y^{[1]}=y'-h(x)\sqrt{\rho}y ,   \qquad \text{  где }h=\int \frac{q}{\sqrt{\rho}} \, dx .$$ А именно,
	 \begin{equation}\label{7}
	y^{[1]}_{\pm}(x, \lambda)= \pm i\lambda \left[\rho(x)\right]^{1/4}\exp\left( \frac 12 P(x) \pm i \lambda\, \int_a^x\sqrt{\rho(\tau)}\, d\tau\right)
\left(1+ \psi_\pm(x, \lambda)\right),
	 \end{equation}
где функции $\psi_\pm$  обладают таким же свойством, как функции $\varphi_\pm$. Утверждение теоремы сохраняется, если вместо
полуплоскости $\im\lambda \geqslant -r$ рассматривать полуплоскость $\im\lambda \leqslant r$.
}

{\sc Доказательство.}
{\sl Шаг 1.} Проведем стандартную замену $$t=\int\limits_a^x \sqrt{\rho (\xi)} \, dx  .$$
	Тогда
$$
 y'_x=y'_t  t'_x=y'_t  \sqrt{\rho(x)},
$$
$$
y''_{xx}=(y'_t   \sqrt{\rho(x)} )'_x=(y'_t)'_x  \sqrt{\rho(x)}+y'_t\frac{\rho '_x}{2\sqrt{\rho (x)}}=y''_{tt}
  \rho+\frac{1}{2}y'_t \rho '_t.
$$
 Подставим эти выражения в уравнение \eqref{3} и разделим на $\rho$. Получим
	  \begin{equation} \label{8}
	  -y''+\left(\frac{p}{\sqrt{\rho}}-\frac{1}{2}\frac{\rho'}{\rho}\right)y'+\frac{q}{\rho}y=\lambda^2 y.
	  \end{equation}
	  Здесь все производные берутся по переменной $t \in [0,h]$, где $h=\int\limits_a^b \sqrt{\rho (\xi)}d\xi$.
	 Заметим, что по условию функции
	 $p$ и $ \rho' $ от переменной $x$ принадлежат пространствам $L_2[a,b]$ и $L_1[a,b]$. В силу равенства $dt=\sqrt{\rho(x)} dx$ эти свойства сохраняются, если их рассматривать как функции от переменной   $t \in [0,h]$. Отметим здесь же, что из абсолютной непрерывности функции $\rho$ и условия $ \rho(x)>0$  при  $x \in [a,b]$ следует ограниченность функции $\rho$  и ее равномерная отделенность от нуля.
	
	{\sl Шаг 2. } Положим
	$$\sigma (t)= \int \frac{q(t)}{\rho (t)}\, dt, $$
	где первообразная берется по переменной $t$ в смысле теории распределений. Заметим, что $\sigma \in L_2[0,h]$. Действительно, по условию $u=\int q(x) \, dx \in L_2[a,b]$. Тогда
$$
\left( \frac{u}{\sqrt{\rho} } \right)'_t=\left( \frac{u}{\sqrt{\rho} } \right)'_x  \frac{1}{\sqrt{\rho} }=\left( \frac{u'_x}{\sqrt{\rho}} - \frac{\rho'_x u}{2 \rho ^{3/2} } \right) \frac{1}{\sqrt{\rho} }=  \frac{q}{\rho } - \frac{\rho'_x u}{2 \rho ^2 }.
$$
Функция $u/{\sqrt{\rho}} \in L_2[0,h]$, поэтому левая часть последнего равенства принадлежит $ W_2^{-1}[0,h].$ Так как по условию  $\rho'_x u/{2 \rho ^2} \in L_1[0,h] \subset  W_2^{-1}[0,h],$ то $q/{\rho} \in W_2^{-1}[0,h]$, а тогда $\sigma \in L_2[0,h]$.		 

Далее воспользуемся приемом из работы \cite{8}. Введем квазипроизводную
\begin{equation}\label{9}
y^{[1]}=y'-\sigma y
\end{equation}
и перепишем уравнение \eqref{8} в виде
\begin{equation}\label{10}
-(y^{[1]})'+\left(\frac{p}{\sqrt{\rho}}-\frac{1}{2}\frac{\rho'}{\rho}-\sigma\right) y^{[1]}+\left(-\sigma^2+\left(\frac{p}{\sqrt{\rho}}-\frac{1}{2}\frac{\rho'}{\rho}\right)\sigma\right)y=\lambda^2y.
\end{equation}
Обозначим
$$
f= \frac{p}{\sqrt{\rho}}-\frac{1}{2}\frac{\rho'}{\rho}, \qquad g=\left(\frac{p}{\sqrt{\rho}}-\frac{1}{2}\frac{\rho'}{\rho}\right)\sigma -\sigma^2.
$$	
Отметим, что  в силу условий \eqref{4}, \eqref{5} обе функции $f$ и $g$ принадлежат пространству $L_1[0,h]$. Уравнения \eqref{9} и \eqref{10} эквивалентны системе уравнений
\begin{equation}\label{11}
\begin{pmatrix} y\\ y^{[1]} \end{pmatrix}'
 = \Lambda\begin{pmatrix} y\\ y^{[1]} \end{pmatrix},
\quad\text{где} \ \ \Lambda = \begin{pmatrix} \sigma & 1 \\ -\lambda^2 +g & f-\sigma \end{pmatrix}.
\end{equation}

Известно (см., например \cite[\S 18]{5}), что такая система уравнений (при условии, что все элементы матрицы суммируемые функции) имеет два линейно независимых решения, причем обе компоненты $y$ и $y^{[1]}$ этих решений абсолютно непрерывны. Далее под решениями исходного уравнения \eqref{3} мы будем понимать функции $y=y(t(x))$, где $y(t)$ -- решение задачи \eqref{11}.

{\sl Шаг 3.} Положим $y=y_1, \, y^{[1]}=y_2$ и перепишем	систему \eqref{11}  в виде	
				\begin{equation}\label{12}
		\begin{pmatrix}	y_1 \\
			y_2\end{pmatrix}'
					=\Lambda_0\begin{pmatrix}
			y_1\\
			y_2 \end{pmatrix}
			+\Lambda_1\begin{pmatrix}
			y_1\\ y_2 \end{pmatrix}, \quad \text{где} \quad \Lambda_0= \begin{pmatrix}
		0&1 \\ -\lambda^2 & 0 \end{pmatrix},
\ \, \Lambda_1= \begin{pmatrix}
			\sigma&0\\
		g& f-\sigma \end{pmatrix}.
		\end{equation}
		
		Решение системы $\mathbf y=\Lambda_0 \mathbf y$ находим явно
		
		$$\mathbf{y}:= \begin{pmatrix}	y_1 \\
			y_2\end{pmatrix} =  M \begin{pmatrix}
			C_1 \\ C_2		
			\end{pmatrix},
			\qquad M= \begin{pmatrix}
			e^{\mu t} & e^{-\mu t}\\
			\mu e^{\mu t} & -\mu e^{-\mu t}\end{pmatrix}, \quad \mu= -i\lambda,
$$
где $C_1, C_2$ -- произвольные постоянные.
Далее воспользуемся методом вариации постоянных. Найдем функции  $C_1=C_1(t),\  C_2=C_2(t)$, которые являются решениями системы
$$
M	\begin{pmatrix}
			C_1\\
			C_2
	\end{pmatrix}' =\Lambda_1
			\begin{pmatrix}
			y_1\\
			y_2
	\end{pmatrix} =
		\begin{pmatrix}
		\sigma & 0\\
		g & f-\sigma
	\end{pmatrix}		
		\begin{pmatrix}
		y_1\\
		y_2
		\end{pmatrix}, \quad M^{-1}=\frac{1}{2\mu}
		\begin{pmatrix}
		\mu e^{-\mu t} & e^{-\mu t}  \\
		\mu e^{\mu t} & -e^{\mu t}
 		\end{pmatrix}.
 $$
Тогда
$$		\begin{pmatrix}
			C_1\\
			C_2
		\end{pmatrix}'
		= M^{-1}\Lambda_1
		\begin{pmatrix}
		y_1\\
		y_2
		\end{pmatrix}
		=K
		\begin{pmatrix}
		y_1\\
		y_2
		\end{pmatrix}, \quad K=\frac{1}{2}	\begin{pmatrix}
		 e^{-\mu t}(\sigma+\mu^{-1}g)& \mu^{-1}e^{-\mu t} (f-\sigma) \\
		 e^{\mu t}(\sigma-\mu^{-1}g) & -\mu^{-1}e^{\mu t} (f-\sigma)
		\end{pmatrix}.
$$
Следовательно, решение системы имеет вид
\begin{equation}\label{13}
	\begin{pmatrix}
	y_1\\
	y_2
	\end{pmatrix}
	=M(t)
	\begin{pmatrix}
	C_1(t)\\
	C_2(t)\end{pmatrix}
	=M(t)
	\begin{pmatrix}	C_1^0\\
	C_2^0
	\end{pmatrix}
	+M(t)\int\limits_0^t K(\xi)
	\begin{pmatrix}
	y_1\\
	y_2
	\end{pmatrix}
	 \, d\xi,
\end{equation}
где $C_1^0, \ C_2^0$ -- произвольные постоянные.

		Положим $C_1^0=1, C_2^0=0$ и сделаем замену $y_1 = e^\mu z_1, y_2 = \mu e^\mu z_2.$
Тогда уравнение \eqref{13}  запишется в виде
$$	
\begin{pmatrix}
				z_1\\
				z_2
			\end{pmatrix}
			=M_0^{-1}(t)M(t)
			\begin{pmatrix}
			1\\
				0
			\end{pmatrix}
			+M_0^{-1}(t)M(t)\int\limits_0^t K(\xi) M_0(\xi)
			\begin{pmatrix}
				z_1(\xi)\\
				z_2(\xi)
             \end{pmatrix}
			 \, d\xi.
$$
Здесь
$$
             M_0=
			\begin{pmatrix}
			e^{\mu t}& 0 \\
		0 & \mu e^{\mu t}
			\end{pmatrix}, \quad M_0^{-1}M=
			\begin{pmatrix}
			1& e^{-2 \mu t} \\
			1 & -e^{-2 \mu t}
			\end{pmatrix},
$$
$$
			  KM_0=\frac{1}{2}
			\begin{pmatrix}
		(\sigma+\mu^{-1}g)&  (f-\sigma) \\
			e^{2\mu \xi}(\sigma-\mu^{-1}g) & -e^{2\mu \xi} (f-\sigma)
			\end{pmatrix}.
$$
Теперь последнее уравнение можно записать в виде
		\begin{equation} \label{14}
				\begin{pmatrix}
				z_1\\
				z_2
				\end{pmatrix}=
				\begin{pmatrix}
				1\\
				1
				\end{pmatrix}
				+\mathcal A\begin{pmatrix}
				z_1\\
				z_2
				\end{pmatrix}
				+\mathcal B\begin{pmatrix}
				z_1\\
				z_2
				\end{pmatrix},
	\end{equation}
	где
\begin{equation}\label{A}
         \mathcal A\begin{pmatrix}
			z_1\\
			z_2	\end{pmatrix}
			=\frac{1}{2}\int\limits_0^t
			\begin{pmatrix}
			\sigma+\mu^{-1}g&  f-\sigma \\
			\sigma+\mu^{-1}g & f-\sigma
			\end{pmatrix}
			\begin{pmatrix}
			z_1(\xi)\\
			z_2(\xi)
			\end{pmatrix}
		     \, d\xi,
\end{equation}
\begin{equation}\label{B}
           \mathcal B \begin{pmatrix}
			z_1\\
			z_2
			\end{pmatrix}
			=\frac{1}{2}\int\limits_0^t
			\begin{pmatrix}
			e^{-2\mu (t-\xi)}(\sigma-\mu^{-1}g) &  -e^{-2\mu (t-\xi)}(f-\sigma) \\
			-e^{-2\mu (t-\xi)}(	\sigma-\mu^{-1}g) & 	e^{-2\mu (t-\xi)}(f-\sigma )
			\end{pmatrix}
			\begin{pmatrix}
			z_1(\xi)\\
			z_2(\xi)
			\end{pmatrix}\, d\xi.
\end{equation}
{\sl Шаг 4.} Далее нам потребуются вспомогательные предложения. Разные варианты этих предложений известны, но для нужных нам версий мы не нашли подходящей ссылки, поэтому для полноты изложения приведем доказательства.
		
{\sc Лемма 1} {\it Пусть $B[0,h]$ -- пространство ограниченных функций с нормой $\|z(x)\|=\sup\nolimits_{x \in [0,h]} |z(x)|$. Пусть функция $K(x, \xi)$  при всех $x\in[0,h]$ суммируема по $\xi$ на отрезке  [0,h] и подчинена оценке
$$
|K(x, \xi)| \leqslant C f(\xi), \quad \text{где} \ f \in L_1 [0,h]
$$
с постоянной $C$, не зависящей от $x$  и $\xi$. Тогда оператор
$$
(Tz)(x)=\int\limits_0^x K(x, \xi) z(\xi) \, d\xi
$$
ограничен в пространстве $B[0,h]$  и его спектральный радиус равен нулю.	}
		
{\sc Доказательство.} Ограниченность оператора вытекает из оценки
$$
|(Tz)(x)|\leqslant C F (x)\, \|z(x)\|, \quad \text{где} \ F (x)=\int\limits_0^x f(\xi) \, d\xi.
$$
Далее используем метод индукции. Предположим, что уже доказана оценка
$$
|(T^{n-1}z)(x)|\leqslant \frac{C^{n-1}}{(n-1)!} F^{n-1} (x)\, \|z(x)\|.
$$
Тогда
\begin{multline*}
|(T^n z)(x)|\leqslant \frac{C^{n}}{(n-1)!}\int\limits_0^x f(\xi) F ^{n-1} (\xi) \, d\xi\,  \|z(x)\|=\\ \frac{C^{n}}{(n-1)!}\int\limits_0^x  F^{n-1} (\xi) \, dF(\xi) \|z(x)\|=
\frac{C^{n}}{n!} F^n(x)\, \|z(x)\|.
\end{multline*}
Следовательно, $\|T^n\|\leqslant C^nF ^n(h)/{n!}$. Так как спектральный радиус оператора определяется равенством $r(T)=\limsup \limits_{n\to \infty} \|T^n\|^{1/n}$, то получаем $r(T)=0$.
		
		{\sc Лемма 2.} {\it Пусть
$$
\mathcal{T}=
			\begin{pmatrix}
		T_1 &  	T_2 \\
			T_3 & 	T_4
\end{pmatrix}
$$ -- оператор в пространстве $B[0,h] \oplus B[0,h]$, где 	
$$
(T_jz)(x)=\int\limits_0^x K_j(x, \xi) z(\xi) \, d\xi ,
$$
причем ядра $K_j$ подчинены условию Леммы 1, то есть,
$$
|K_j(x, \xi)| \leqslant C_j f_j(\xi), \quad \text{где} \ f_j \in L_1 [0;h], \quad j=1,2,3,4.
$$
Тогда оператор $\mathcal{T}$ ограничен и его спектральный радиус равен 0. Более того,
$$
\|\mathcal{T}^n\| \leqslant\frac{2 (2C)^n F^n(h)}{n!}, \quad \text{где }\  F (x)=\int\limits_0^x f(\xi)\, d\xi,\ f(x)=\max\limits_{1 \leqslant j \leqslant4} f_j(x), \ C=\max\limits_{1\leqslant j \leqslant4} C_j.
$$
Здесь норму вектора
$ z=
			\begin{pmatrix}
			z_1\\
			z_2
			\end{pmatrix}
$  в пространстве $B[0,h] \oplus B[0,h]$ мы определяем равенством $\|z\|=\max (\|z_1\|, \|z_2\|).$ }

{\sc Доказательство.} Заметим, что
			$$\mathcal{T}^n=
				\begin{pmatrix}
				U_1&  	U_2 \\
				U_3 & 	U_4
				\end{pmatrix}
			 $$
где каждый из операторов $U_j$, $j=1,2,3,4,$ есть сумма $2^{n-1} $ слагаемых, а эти слагаемые есть произведения вида $T_{m_1} T_{m_2} \cdots T_{m_n}$ с индексами ${1 \leqslant m_k \leqslant 4}$ (конечно, среди индексов могут быть одинаковые). Тогда
\begin{multline*}
 |(T_{m_1} T_{m_2} \cdots T_{m_n}z)(x)| \\
 \leqslant \int\limits_0^x |K_{m_1}(x,\xi) \int\limits_0^\xi K_{m_2}(\xi,t) \cdots  \int\limits_0^\tau K_{m_n}(\tau,s)z(s) \, ds\, d\tau \cdots dt | \, d\xi\, \|z\|  \\
\leqslant C_{m_1} C_{m_2}\cdots C_{m_n} \int\limits_0^x  f_{m_1} (\xi) \int\limits_0^{\xi}  f_{m_2} (t) \cdots  \int\limits_0^\tau  f_{m_n} (s) \, ds \, d\tau \cdots dt \, d\xi\, \|z\|  \\
\leqslant C^n  \int\limits_0^x  f (\xi) \int\limits_0^{\xi}  f (t) \cdots  \int\limits_0^\tau  f (s) \, ds \, d\tau \cdots dt \, d\xi \|z\| \leqslant \frac{C^n F^n(h)}{n!}\, \|z\|.
\end{multline*}
Очевидно, норма оператора $\mathcal{T}^n$ допускает оценку $\|\mathcal{T}^n\| \leqslant \|U_1\|+ \|U_2\|+ \|U_3\|+ \|U_4\|$.
Но тогда  $\|\mathcal{T}^n\| \leqslant 4 \cdot 2^{n-1} C^n F^n (h)/{n!}$, что доказывает лемму.
			
	Заметим, что если ядра интегральных операторов непрерывны по $x$ как функции со значениями в пространстве $L_1[0,h]$, то соответствующие операторы корректно определены в пространстве непрерывных функций $C[0,h]$. Оценки для норм в этом пространстве, конечно, сохраняются.

{\sl Шаг 5.} Далее будем работать в полуплоскости $\re \mu \geqslant -r$,  где получим
 асимптотические представления  для решений.  Вернемся к представлению \eqref{14} и заметим, что интегральные операторы $\mathcal A$ и $\mathcal B$, определенные интегральными формулами \eqref{A}  и \eqref{B}, таковы, что элементы их матричных ядер допускают  оценку
$$
\leqslant e^{2rh} v(\xi), \qquad \text{где} \  v(\xi)=|f(\xi)|+|g(\xi)|+|\sigma(\xi)|,
$$
равномерно по $0\leqslant \xi\leqslant x\leqslant h$ при $\re \mu \geqslant -r, \, |\mu | \geqslant 1.$
		
		Решение уравнения \eqref{14} запишем в виде
	\begin{equation}\label{15}
			\begin{pmatrix}
			z_1\\
			z_2
			\end{pmatrix}
			=(1-\mathcal A)^{-1}
			\begin{pmatrix}
			1\\
			1
			\end{pmatrix} +\mathcal T
			\begin{pmatrix}
			z_1\\
			z_2
			\end{pmatrix}, \quad \mathcal T=(1-\mathcal A)^{-1} \mathcal B.
\end{equation}

		Заметим, что интегральный оператор $\mathcal A+\zeta\mathcal B$
вольтерров, т.е. его спектральный радиус равен нулю при любом $\zeta \in \mathbb{C}$. Поэтому оператор
$$
(1-\mathcal A-\zeta \mathcal B)=(1-\mathcal A)(1-\zeta\mathcal T)
$$
обратим при любом $\zeta \in \mathbb{C}$. Но тогда $(1- \zeta \mathcal T)$ обратим при любом $\zeta \in \mathbb{C}$, т.е. спектральный радиус оператора $\mathcal T$ равен нулю. Следовательно, при любом $\varepsilon \in (0,1)$ имеем $\|\mathcal T^n\| < \varepsilon^n$, если $n \geqslant N=N(\varepsilon)$. Решение уравнения \eqref{15} теперь запишем в виде
		\begin{equation} \label{16}
			\begin{pmatrix}
			z_1\\
			z_2
			\end{pmatrix}
			=\begin{pmatrix}
			\psi_1\\
			\psi_2
			\end{pmatrix}
			+\sum\limits_{k=1}^{N-1}\mathcal T^k
			\begin{pmatrix}
			\psi_1\\
			\psi_2
			\end{pmatrix}
			+\mathcal T^N(1-\mathcal T)^{-1}
			\begin{pmatrix}
			\psi_1\\
			\psi_2
			\end{pmatrix}, \quad
			\begin{pmatrix}
			\psi_1\\
			\psi_2
			\end{pmatrix}
			=(1-\mathcal A)^{-1}
			\begin{pmatrix}
			1\\
			1
			\end{pmatrix}
			\end{equation}
Согласно \eqref{A} имеем
$$      \mathcal A\begin{pmatrix}
		1\\
		1
		\end{pmatrix}
		=\frac{1}{2}\int\limits_0^t
		\begin{pmatrix}
		\sigma+\mu^{-1}g&  f-\sigma \\
		\sigma+\mu^{-1}g & f-\sigma
		\end{pmatrix}
		\begin{pmatrix}
		1\\
		1 \end{pmatrix}
		\, d\xi=
		\frac 12\begin{pmatrix}
		F(t)\\
		F(t)
		\end{pmatrix}
		+O\left(\mu^{-1}\right), \quad F(t)=\int\limits_0^t f(\xi) \, d\xi.
$$
	Тогда
\begin{equation}\label{main}
        \mathcal A^n
		\begin{pmatrix}
		1\\
		1
		\end{pmatrix}
		=\frac{1}{2^n n!}
		\begin{pmatrix}
		F^n(t)\\
		F^n(t)
		\end{pmatrix} + O(\mu^{-1}) \frac{C^{n}}{n!}, \quad (1-\mathcal A)^{-1}
		\begin{pmatrix}
		1\\
		1
		\end{pmatrix} =
		\begin{pmatrix}
		e^{\frac 12 F(t)}\\
		e^{\frac 12 F(t)}
		\end{pmatrix} +O(\mu^{-1}).
\end{equation}
Поясним, как получается эти соотношения. Легко видеть, что
$$
\mathcal A^n
		\begin{pmatrix}
		1\\
		1
		\end{pmatrix} =
\frac 1{2^n} \begin{pmatrix} (V_1 + \mu^{-1} V_2)^n \\ (V_1+ \mu^{-1} V_2)^n\end{pmatrix}
\begin{pmatrix} 1\\1\end{pmatrix},$$ где $$ V_1 z = \int\limits_0^t f(\xi) z(\xi)\, d\xi,
\ \ V_2 z =   \int\limits_0^t g(\xi) z(\xi)\, d\xi.
$$
Имеем $$(V_1+\mu^{-1} V_2)^n = V_1^n + \mu^{-1} U_n,$$  где $U_n$  есть
есть сумма $2^n-1$  слагаемых  вида $ \mu^{-j} V_{k_1}V_{k_2} \cdots V_{k_n}$.  Здесь $j\geqslant 0$,  а
индексы $k_j$  принимают значения 1 и 2. Повторяя рассуждения Леммы 2, получаем при $|\mu|\geqslant 1$  оценки
$\|U_n\| \leqslant C^n/ n!$. Из этих оценок  следуют соотношения \eqref{main}.

Согласно  введенным обозначениям получаем
	\begin{multline}\label{17}
	\frac 12 F(\xi)=\frac{1}{2}\int\limits_0^{\xi} \frac{p(x(t))}{\sqrt{\rho (x(t))}}\, dt-\frac{1}{4}\int\limits_0^{\xi} \,
 \frac{\rho'_t}\rho\, dt \\ =\frac{1}{2}\int\limits_0^{\xi} \frac{p(x(t))}{\sqrt{\rho (x(t))}}\, t'_x \,dx(t)
 -\frac{1}{4}\int\limits_0^\xi\, \frac{\rho'_x}{t'_x\rho} \, t'_x dx=
\frac{1}{2} P(\xi)+\ln \frac {1}{\sqrt[4]{\rho(\xi)} }.
	\end{multline}

	Заметим, что для любой фиксированной функции $\psi (\xi) \in L_1(0,h)$ выполняется соотношение
	\begin{equation} \label{18}
\left|\int\nolimits_0^t e^{2\mu (t-\xi)} \psi (\xi) \, d\xi\right| \rightarrow 0 \quad \text{при} \quad|\mu| \rightarrow \infty
	\end{equation}
равномерно по $t \in [0,h]$ в полуплоскости $\re \mu\geqslant -r$. Доказательство этого факта стандартно. Сначала с помощью интегрирования по частям получаем, что для гладкой функции $\psi$ этот интеграл есть $O(\mu^{-1})$ при $|\mu| \rightarrow \infty$, а тогда утверждение \eqref{18} получается из того, что суммируемую функцию можно с произвольной точностью в метрике пространства $L_1$ приблизить функцией гладкой.

	Обратимся к представлению \eqref{16}. Фиксируем произвольное число $\varepsilon$ и выберем число $N$ так, чтобы
$$
\|\mathcal T^N(1-\mathcal T)^{-1}(1-\mathcal A)^{-1}\|<\varepsilon .
$$
Тогда норма в пространстве $C[0,h]\oplus C[0,h]$ последнего слагаемого в \eqref{16} оценивается числом $\varepsilon$. В силу
\eqref{B} и \eqref{18} имеем (здесь $\mathcal B=\mathcal B(\mu)$)
	$$\left\|\mathcal B \begin{pmatrix}
		\psi_1\\
		\psi_2
		\end{pmatrix}
		\right\| \rightarrow 0 \quad \text{при} \  |\mu| \rightarrow \infty \quad \text{в полуплоскости}\quad \re \mu\geqslant -r.
$$
Но тогда при любом фиксированном $k\geqslant 1$ получаем
$$
        \left\|\mathcal T^k
		\begin{pmatrix}
		\psi_1\\
		\psi_2
		\end{pmatrix}\right\|=
        o(1) \quad \text{при} \ |\mu| \rightarrow \infty , \quad  \re \mu\geqslant -r.
$$
	Следовательно нормы всех слагаемых, начиная со второго, в правой части \eqref{16} оцениваются числом $\varepsilon/N$ при условии, что $|\mu |>r_0$, а $r_0$ достаточно большое число. Тогда из \eqref{16} и \eqref{17} получаем
\begin{equation}\label{Ok}
		\begin{pmatrix}
		z_1\\
	    z_2
		\end{pmatrix}=
		\begin{pmatrix}
		e^{\frac 12 F(t)}\\
		e^{\frac 12 F(t)}
		\end{pmatrix}
		(1+o(1)) \quad \text{при} \ |\mu| \rightarrow \infty , \quad  \re \mu\geqslant -r,
\end{equation}
где функция $F$ определена  в \eqref{17}.
Теперь утверждение теоремы для решения $y_+(x,\lambda)$ получается из равенств
$$
		\begin{pmatrix}
		y(x,\lambda)\\
		y^{[1]}(x,\lambda)
		\end{pmatrix}=
		\begin{pmatrix}
		e^{\mu t}& 0 \\
		0 & \mu e^{\mu t}
		\end{pmatrix}
		\begin{pmatrix}
		z_1\\
		z_2
		\end{pmatrix}, \qquad \lambda = i\mu.
$$

{\sl Шаг 6.} Докажем существование решения $y_-(x,\lambda)$  в полуплоскости $\im \lambda \geqslant -r$. Систему \eqref{11} мы представили ранее в виде интегрального уравнения \eqref{13}. Заменив в \eqref{13} постоянные $C_1^0, C_2^0$ мы запишем интегральное уравнение в другой форме
$$
		\begin{pmatrix}
			y_1\\
			y_2
		\end{pmatrix}=
         M(t)\begin{pmatrix}
			C_1^0\\
			C_2^0
		\end{pmatrix} - M(t)\int\limits_t^h \, K(\xi)
		\begin{pmatrix}
			y_1 (\xi)\\
			y_2 (\xi)
		\end{pmatrix} \, d\xi.
$$
Здесь  $C_1^0, C_2^0$ снова произвольные постоянные, но уже другие, нежели в \eqref{13}.
Положим $C_1^0=0$,  $C_2^0=1$  и  сделаем замену $y_1=e^{-\mu t}z_1, \ y_2= - \mu e^{-\mu t}z_2. $ Тогда уравнение перепишется в виде
\begin{equation}\label{last}
		\begin{pmatrix}
		z_1\\
		z_2
		\end{pmatrix}
		=M_1^{-1}(t)M(t)
		\begin{pmatrix}
		0\\
		1
		\end{pmatrix}
		- M_1^{-1}(t)M(t)\int\limits_t^h \, K(\xi) M_1(\xi)
		\begin{pmatrix}
		z_1(\xi)\\
		z_2(\xi)
		\end{pmatrix}\, d\xi,
\end{equation}
Здесь
$$           M= \begin{pmatrix}
			e^{\mu t} & e^{-\mu t}\\
			\mu e^{\mu t} & -\mu e^{-\mu t}\end{pmatrix},
            \ \  M_1=
			\begin{pmatrix}
			e^{-\mu t}& 0 \\
		0 & -\mu e^{-\mu t}
			\end{pmatrix}, \quad M_1^{-1}M=
			\begin{pmatrix}
			 e^{2 \mu t}& 1 \\
			-e^{2 \mu t}& 1
			\end{pmatrix},
$$
$$
			  K M_1=\frac{1}{2}
			\begin{pmatrix}
             e^{-2\mu \xi}(\sigma+\mu^{-1}g) & -e^{-2\mu \xi} (f-\sigma)\\
		     \sigma-\mu^{-1}g &  f-\sigma
			\end{pmatrix}.
$$
Теперь уравнение \eqref{last}  можно представить в виде
$$	
		\begin{pmatrix}
		z_1\\
		z_2
		\end{pmatrix}=	
		\begin{pmatrix}
		1\\
		1
		\end{pmatrix}
		+\mathcal A_1\begin{pmatrix}
		z_1\\
		z_2
		\end{pmatrix}
		+ \mathcal B_1 \begin{pmatrix}
		z_1\\
		z_2
		\end{pmatrix},
$$
где
$$          \mathcal A_1\begin{pmatrix}
			z_1\\
			z_2
			\end{pmatrix}=
			-\frac{1}{2}\int\limits_t^h
			\begin{pmatrix}
			\sigma-\mu^{-1}g &  f-\sigma \\
			\sigma-\mu^{-1}g & f-\sigma
			\end{pmatrix}
			\begin{pmatrix}
			z_1(\xi)\\
			z_2(\xi)
			\end{pmatrix} \, d\xi,
$$
$$
             \mathcal B_1
			\begin{pmatrix}
			z_1\\
			z_2
			\end{pmatrix}=-\frac{1}{2}\int\limits_t^h
			\begin{pmatrix}
			e^{2\mu (t-\xi)}(\sigma+\mu^{-1}g)&  	-e^{2\mu (t-\xi)}(f-\sigma) \\
			-e^{2\mu (t-\xi)}(	\sigma+\mu^{-1}g) & 	e^{2\mu (t-\xi)}(f-\sigma )
			\end{pmatrix}
			\begin{pmatrix}
			z_1(\xi)\\
			z_2(\xi)
			\end{pmatrix} \, d\xi.
$$
Очевидно, операторы $\mathcal A_1$ и  $\mathcal B_1$ также вольтерровы в пространстве $C[0,h]\oplus C[0,h]$. Кроме того,
$$
|e^{2\mu (t-\xi)}| \leqslant e^{2rh} \quad  \text{при}\ \,   \re \mu \geqslant-r \ \, \text{и}\ \, \xi \in [t,h], \ t\in [0,h].
$$
Поэтому все прежние рассуждения сохраняются  и мы получаем, что решение уравнения \eqref{last} представимо в виде
\eqref{Ok}. Тогда утверждение теоремы для решения $y_-(x,\lambda)$  получается из равенств
 $$
		\begin{pmatrix}
		y_-(x,\lambda)\\
		y_-^{[1]}(x,\lambda)
		\end{pmatrix}=
		\begin{pmatrix}
		e^{-\mu t}& 0 \\
		0 & -\mu e^{-\mu t}
		\end{pmatrix}
		\begin{pmatrix}
		z_1\\
		z_2
		\end{pmatrix}, \qquad \lambda=i\mu.
$$
Асимптотические представления решений  в полуплоскости $\im \lambda \leqslant r$ получаются
аналогично. Этим завершается доказательство теоремы.


\end{document}


\author{  V.\,E.~Vladykina \thanks{Lomonosov Moscow State University, e-mail: vika-chan@mail.ru}
}
 \title{Asymptotics of fundamental solutions of the Sturm-Liouville equation for the large spectral parameter}
\date{16.05.2017}

\maketitle

\footnotetext{The work is supported by the grant of President of RF, № NSh- 7461.2016.1.}

{\bf Abstract}.  We consider the Sturm-Liouvile equation
\begin{equation*}
-(r^2y')'+py'+qy=\lambda ^2 \rho^2 y, \qquad x\in [a,b] \subset  \mathbb{R},
\end{equation*}
 where $\lambda^2$ is the large parameter,
  $r$ and $\rho$  are positive functions, while  $p$ and $q$ are complex valued  ones. It is assumed that
\begin{equation*}
 p\in L_1[a,b],\quad q\in W_2^{-1}[a,b], \quad  \rho ,r \in AC[a,b] =W_1^1[a,b],
\end{equation*}
moreover,
\begin{equation*}
\rho'u, r'u, pu \in L_1[a,b], \quad \text{where \,} u=\int q \, dx,
\end{equation*}
and the antiderivative is understood in the sense of distributions.
We find the main term of two fundamental solutions of this equation
with the large parameter $\lambda$ in the upper and the lower complex $\lambda$-half-planes.
The main matter of this note is that the asymptotic formulae are obtained under minimal smoothness assumptions
on the involved functions, moreover,  the function $q$ is allowed to be a distribution of the first order singularity.
\medskip


\bigskip

This note is the supplement to our paper \cite{ShVl}. Further we always use the notations of this paper.
The historical background can be also found in \cite{ShVl}.

  We deal with the Sturm-Liouville equation
\begin{equation} \label{0}
-(r^2y')'+py'+qy=\lambda ^2 \rho^2 y, \qquad x\in [a,b] \subset  \mathbb{R},
\end{equation}
where $\lambda$ is the large parameter,  $r$ and $\rho$  are positive functions, while $p$ and $q$
are complex valued ones. We assume that
\begin{equation}
\label{usl}
 p\in L_1[a,b],\quad q\in W_2^{-1}[a,b], \quad  \rho ,r \in AC[a,b]=W_1^1[a,b],
\end{equation}
moreover,
\begin{equation} \label{dopusl}
\rho'u, r'u, pu \in L_1[a,b], \quad \text{where \,} u=\int q \, dx
\end{equation}
and the antiderivative is understood in the sense of distributions.

\textbf{Main theorem}. {\sl Let the coefficients of equation \eqref{0} are subject conditions
\eqref{usl} and \eqref{dopusl}. Then there exist fundamental solutions $y_\pm$ of equation \eqref{0}
having a representation of the form
\begin{equation} \label{6_1}
y_{\pm}(x, \lambda)= \frac{1}{\sqrt{r(x)\rho(x)}}\exp\left( \frac 12 \int_a^x \frac{p(\tau)}{r^2(\tau)}\, d\tau \pm i \lambda\, \int_a^x\frac{\rho(\tau)}{r(\tau)}\, d\tau \right)
\left(1+ \varphi_\pm(x, \lambda)\right).
\end{equation}
Here the functions $\varphi_{\pm}$ are analytic in a half-plane $\Pi^+_s = \{\lambda\in \mathbb C\vert\  \im\lambda \geqslant -s\}$ provided that $|\lambda|>R$ and
 $$
  |\varphi _+\left(x, \lambda\right)|+ |\varphi _-(x, \lambda)| = o(1) \quad \text{as} \ \, |\lambda | \to \infty \ \,
   \lambda \in \Pi^+_s,
  $$
   uniformly with respect to $x\in [a,b]$.Here positive  number $s$ can be taken arbitrary large. Asymptotic formulae  \eqref{6_1} can be differentiated term-by-term if the derivative is replaced by the quasiderivative	
   $$
   y^{[1]}=y'-h(x)\frac{\rho(x)}{r(x)}y,\quad \text{where}\, h(x)=\int \frac{q}{r \rho} \, dx.
   $$
   Namely,
	 \begin{equation}\label{7_1}
	y^{[1]}_{\pm}(x, \lambda)= \pm i\lambda\sqrt{\frac{\rho}{r^3}}\exp\left( \frac 12 \int_a^x \frac{p}{r^2} \pm i \lambda\, \int_a^x\frac{\rho}{r}\right)
\left(1+ \psi_\pm(x, \lambda)\right),
\end{equation}
where the functions  $\psi_\pm$  possess the same property as the functions $\varphi_\pm$. The assertion of the theorem is preserved if, instead of the half-plane $\Pi^+_s$, the half-plane ${\Pi^-_s=  \{\lambda\in \mathbb C\vert\  \im\lambda \leqslant s\}}$ is considered.
}

In the case $r\equiv 1, p \in L_2[a,b]$ this theorem is proved in the paper \cite{ShVl}. The objective of this note is to prove this result under more general assumptions.

{\bf Proof}. Let us carry out formal transformations of \eqref{0}. We have
 $$-r^2y''-2rr'y'+py'+qy=\lambda ^2 \rho^2 y,$$
$$-y''+\left(-2\frac{r'}{r}+\frac{p}{r^2}\right)y'+\frac{q}{r^2}=\lambda ^2 \frac{\rho^2}{r^2} y.$$
Using the standard substitution
\begin{equation}\label{subst}
t=\int\limits_a^x \frac{\rho (\xi)}{r(\xi)} \, d\xi  .
\end{equation}
we obtain
\begin{equation} \label{eq_p}
  -y''_{tt}+\left(-\left(\frac{\rho}{r}\right)'_t\cdot\frac{r}{\rho}-2\frac{r'}{r}+\frac{p}{r\rho}\right)y'_t+\frac{q}{\rho^2}y=\lambda^2 y, \qquad
t \in [0,h],\ h=\int\limits_a^b \frac{\rho (\xi)}{r(\xi)} d\xi.
\end{equation}
	
Consider the function
	  	$$\sigma (t)= \int \frac{q(t)}{\rho^2 (t)}\, dt, $$
where the antiderivative is taken with resprct to the variable  $t$ in the sense of distributions.
Note, that  $\sigma \in L_2[0,h]$. Indeed, by assumption \eqref{usl} we have $u=\int q(x) \, dx \in L_2[a,b]$.
Then
\begin{equation} \label{potential}
\left( \frac{u}{r\rho } \right)'_t=\frac{r}{\rho } \left( \frac{u}{r\rho } \right)'_x  =\left( \frac{u'_x}{r\rho} - \frac{ur'_x}{r^2\rho}-\frac{u\rho'_x}{r\rho^2} \right) \frac{r}{\rho} =  \frac{q}{\rho^2 } - \frac{ur'_x}{r\rho^2}-\frac{\rho'_x u}{\rho ^3 }.
\end{equation}

The function $u/{r\rho} \in L_2[a,b]$ since it is the product of a function from $L_2$ and a bounded measurable function . By virtue \eqref{subst} we have also   $u/{r\rho} \in L_2[0,h]$ similarly, therefore the left hand-side of the last equality
  belongs to the space $ W_2^{-1}[0,h].$  By virtue of conditions  \eqref{dopusl}
$$
\rho'_x u, r'_x u \in L_1[0,h], \quad \text{hence} \ \ \frac{ur'_x}{r\rho^2}+\frac{\rho'_x u}{\rho ^3 }\in L_1[0;h]\subset  W_2^{-1}[0,h].
$$
This implies  $q/{\rho^2} \in W_2^{-1}[0,h]$, and then  $\sigma \in L_2[0,h]$.		

Let us introduce the quasiderivative
\begin{equation}\label{eq_1}
y^{[1]}=y'-\sigma y
\end{equation}
and rewrite equation \eqref{eq_p} in the form
\begin{equation}\label{eq_2}
-(y^{[1]})'+\left(-\left(\frac{\rho}{r}\right)'\cdot\frac{r}{\rho}-2\frac{r'}{r}+\frac{p}{r\rho}-\sigma\right) y^{[1]}+\left(-\sigma^2+\left(-\left(\frac{\rho}{r}\right)'\cdot\frac{r}{\rho}-2\frac{r'}{r}+\frac{p}{r\rho}\right)\sigma\right)y=\lambda^2y.
\end{equation}
Denote
$$
f= -\left(\frac{\rho}{r}\right)'\cdot\frac{r}{\rho}-2\frac{r'}{r}+\frac{p}{r\rho}, \qquad g=\left(-\left(\frac{\rho}{r}\right)'\cdot\frac{r}{\rho}-2\frac{r'}{r}+\frac{p}{r\rho}\right)\sigma -\sigma^2.
$$	
Obviously, $f\in L_1[0,h]$. Let us show that $g\in L_1[0,h]$.  This follows from  \eqref{usl},\eqref{dopusl}.
We have to check only that
$$
\frac{r'}{r}\sigma, \left(\frac{\rho}{r}\right)'\cdot\frac{r}{\rho}\sigma \in L_1[0,h].
$$
To prove this statement it is sufficient to show that $\rho'_t \sigma \in L_1[0,h]$ and $r'_t \sigma \in L_1[0,h]$.
Indeed, similarly \eqref{potential} we have
$$\rho'_t \sigma =\rho'_t \int \frac{q}{\rho r} \, dx=\rho'_t \int \left(\frac{u}{\rho r}\right)'_x \, dx +\rho'_t \int  \frac{u(r\rho )'_x}{r\rho}\, dx=\rho'_t \frac{u}{\rho r}+\rho'_t \int  \frac{u(r\rho )'_x}{r\rho}\, dx.$$
By virtue of  \eqref{dopusl} and \eqref{subst} the first summand in the right hand-side  belongs to $L_1[0,h]$, and so do
the second summand, since this is represented as the product of a function from $L_1[0,h]$ and a function
 from $W^1_1[0,h]= AC[0,h]$. Similarly we have $r'_t \sigma \in L_1[0,h]$. Hence  $g\in L_1[0,h]$.

Equations  \eqref{eq_1} and \eqref{eq_2} are equivalent to the system of equations
\begin{equation}\label{sys_p}
\begin{pmatrix} y\\ y^{[1]} \end{pmatrix}'
 = \Lambda\begin{pmatrix} y\\ y^{[1]} \end{pmatrix},
\quad\text{where} \ \ \Lambda = \begin{pmatrix} \sigma & 1 \\ -\lambda^2 +g & f-\sigma \end{pmatrix}.
\end{equation}
Here the coefficients  $f,g \in L_1[0;h],\, \sigma \in L_2[0;h]$.
Repeating word-by-word the arguments from the paper  \cite{ShVl}, applying the method of constant variation
 and making the substitution
 \begin{equation} \label{z}
 \begin{pmatrix}
		y(x,\lambda)\\
		y^{[1]}(x,\lambda)
		\end{pmatrix}=
		\begin{pmatrix}
		e^{\mu t}& 0 \\
		0 & \mu e^{\mu t}
		\end{pmatrix}
		\begin{pmatrix}
		z_1\\
		z_2
		\end{pmatrix}, \qquad \lambda = i\mu.
 \end{equation}

we get the system of equations for the new functions
	\begin{equation*}
				\begin{pmatrix}
				z_1\\
				z_2
				\end{pmatrix}=
				\begin{pmatrix}
				1\\
				1
				\end{pmatrix}
				+\mathcal A\begin{pmatrix}
				z_1\\
				z_2
				\end{pmatrix}
				+\mathcal B\begin{pmatrix}
				z_1\\
				z_2
				\end{pmatrix}.
	\end{equation*}
Here
\begin{equation*}
         \mathcal A\begin{pmatrix}
			z_1\\
			z_2	\end{pmatrix}
			=\frac{1}{2}\int\limits_0^t
			\begin{pmatrix}
			\sigma+\mu^{-1}g&  f-\sigma \\
			\sigma+\mu^{-1}g & f-\sigma
			\end{pmatrix}
			\begin{pmatrix}
			z_1(\xi)\\
			z_2(\xi)
			\end{pmatrix}
		     \, d\xi,
\end{equation*}
\begin{equation*}
           \mathcal B \begin{pmatrix}
			z_1\\
			z_2
			\end{pmatrix}
			=\frac{1}{2}\int\limits_0^t
			\begin{pmatrix}
			e^{-2\mu (t-\xi)}(\sigma-\mu^{-1}g) &  -e^{-2\mu (t-\xi)}(f-\sigma) \\
			-e^{-2\mu (t-\xi)}(	\sigma-\mu^{-1}g) & 	e^{-2\mu (t-\xi)}(f-\sigma )
			\end{pmatrix}
			\begin{pmatrix}
			z_1(\xi)\\
			z_2(\xi)
			\end{pmatrix}\, d\xi.
\end{equation*}

Let us rewrite the equation in the form
\begin{equation*}
			\begin{pmatrix}
			z_1\\
			z_2
			\end{pmatrix}
			=(1-\mathcal A)^{-1}
			\begin{pmatrix}
			1\\
			1
			\end{pmatrix} +\mathcal T
			\begin{pmatrix}
			z_1\\
			z_2
			\end{pmatrix}, \quad \mathcal T=(1-\mathcal A)^{-1} \mathcal B.
\end{equation*}
\begin{equation} \label{14}
			\begin{pmatrix}
			z_1\\
			z_2
			\end{pmatrix}
			=\begin{pmatrix}
			\psi_1\\
			\psi_2
			\end{pmatrix}
			+\sum\limits_{k=1}^{N-1}\mathcal T^k
			\begin{pmatrix}
			\psi_1\\
			\psi_2
			\end{pmatrix}
			+\mathcal T^N(1-\mathcal T)^{-1}
			\begin{pmatrix}
			\psi_1\\
			\psi_2
			\end{pmatrix}, \quad
			\begin{pmatrix}
			\psi_1\\
			\psi_2
			\end{pmatrix}
			=(1-\mathcal A)^{-1}
			\begin{pmatrix}
			1\\
			1
			\end{pmatrix}.
			\end{equation}
Here
			\begin{equation}\label{15}
			(1-\mathcal A)^{-1}
		\begin{pmatrix}
		1\\
		1
		\end{pmatrix} =
		\begin{pmatrix}
		e^{\frac 12 F(t)}\\
		e^{\frac 12 F(t)}
		\end{pmatrix} +O(\mu^{-1}).
			\end{equation}
Then on the same way as in \cite{ShVl} we get
\begin{equation}\label{16}
		\begin{pmatrix}
		z_1\\
	    z_2
		\end{pmatrix}=
		\begin{pmatrix}
		e^{\frac 12 F(t)}\\
		e^{\frac 12 F(t)}
		\end{pmatrix}
		(1+o(1)) \quad \text{as} \ |\mu| \rightarrow \infty , \quad  \re \mu\geqslant -r,
\end{equation}
Now, let us return to the original variable			
			\begin{multline*}
	\frac 12 F(\xi)=\frac{1}{2}\int\limits_0^{\xi}\left(- \left(\frac{\rho(x(t))}{r(x(t))}\right)'_t \frac{r(x(t))}{\rho(x(t))}-2\frac{\rho '_t}{\rho}+\frac{p}{r\rho}\right)\, dt \\ =-\frac{1}{2}\ln \frac{\rho}{r}-\ln r+\frac{1}{2}\int\limits_0^{\xi}\frac{p}{r\rho}\cdot t'_x \, dx=
	-\frac{1}{2}\ln \frac{\rho}{r}-\ln r+\frac{1}{2}\int\limits_0^{\xi}\frac{p}{r\rho}\cdot \frac{\rho}{r} \, dx=\\
	=	-\frac{1}{2}\ln \frac{\rho}{r}-\ln r+\frac{1}{2}\int\limits_0^{\xi}\frac{p}{r^2} \, dx.
	\end{multline*}
Now by virtue of \eqref{z} -\eqref{15} we obtain representations \eqref{6_1} and  \eqref{7_1} for the solution $y_+$. Repeating the arguments
we get a similar representation for the solution  $y_-$.